\documentclass[a4paper,10pt]{article}
\usepackage{graphicx}
\usepackage{amsfonts}
\usepackage{amsmath}
\usepackage{amsmath}
\usepackage{stmaryrd}
\usepackage{amssymb}
\usepackage{amssymb}
\usepackage[utf8x]{inputenc}
\usepackage{amsfonts}

\newtheorem{theorem}{Theorem}
\newtheorem{lemma}{Lemma}
\newcommand{\beq}[1]{\begin{equation}\label{#1}}
\newcommand{\eeq}{\end{equation}}
\newcommand{\bref}[1]{   (\ref{#1})  }
\newcommand{\pf}{\noindent{\it Proof}. }

\newtheorem{definition}{Definition}
\newcommand{\m}{\mathbb}
\newcommand{\rem}{\medskip \noindent {\bf Remark: }}
\newcommand{\bea}{\begin{eqnarray*}}
\newcommand{\eea}{\end{eqnarray*}}

\begin{document}

\begin{center}
{\bf\Large From Uniform Continuity to  Absolute Continuity}
\\[3ex]
{\bf Kai Yang, $\,\,\,$ Chenhong Zhu}
\end{center}

The notion of uniform continuity emerged slowly in the lectures of
Dirichlet(1854) and of Weierstrass(1861) \cite{uniformhistory}. Then
in 1905, Vitali established the absolute continuity for a class of
functions in the paper ``Sulle funzioni integrali"
\cite{absolutehistory}. Sooner or later, several equivalent and
sufficient conditions for a function to be absolutely continuous
were derived (see \cite{abs}), which depend on several results of
measure theory and integration. For example, Banach-Zarecki
Criterion\cite{ABS} states that, ``$f$ is absolutely continuous on
$[a,b]$ if and only if $f$ is continuous and of bounded variation on
$[a,b]$ and maps null sets to null sets. '', and in particular,
every Lipschitz continuous function is absolutely continuous.

Absolute continuity implies uniform continuity, but generally not
vice versa. However, under certain conditions (piecewise convexity),
uniform continuity will also imply absolute continuity. In this
short note, we will present a sufficient condition for a uniformly
continuous function to be absolutely continuous.
\begin{definition}{\bf (piecewise convex function)}
A function $f$  defined on an interval $I_{a,b}$ of the real line is
piecewise convex, if there exits a finite partition
$P=\{a_i\}_{i=0}^{N}$ such that on each subinterval $[a_i,a_{i+1}]$
{\small $(i=0,\dots,$} {\small$ N-1)$}, $f$ is concave or convex.
\end{definition}
\noindent{\bf Remark:} Here $I_{a,b}$ is an interval of $\m{R}$ with
$a,b$ as the endpoints and it can be open, closed or half-open;
moreover, $a,b$ can also take the value $\pm\infty$ when $a$ or $b$
is not contained in $I_{a,b}$ (in this case, we use $(a_0,a_{1}]$ or
$[a_{N-1},a_{N})$ instead of the closed subinterval).
\begin{theorem}\label{th1}
 For a uniformly continuous function $f$ defined on $I_{a,b}$ of $\m{R}$, if it is piecewise convex, then it is also absolutely
continuous on $I_{a,b}$.
\end{theorem}
\noindent{\bf Remark:}  We can verify that if $I_{a,b}=[a,b]$, then
the conditions in Theorem \ref{th1} satisfy the Banach-Zarecki
Criterion. However, our attempt is based on some elementary
properties of uniform continuity and convexity. In particular, one
simple example is that $f(x)=\sqrt{x}$ is absolutely continuous (not
Lipschitz continuous) on $[0,c]$. Moreover, the converse statement
of this theorem is false. One proper counterexample is
$f(x)=x^2\sin(1/x)$ on $[0,1]$ (define $f(0)=0$). In addition, the
cantor function\cite{cantor} is uniformly continuous but not
absolutely continuous.
\begin{figure}[h]
\begin{center}
\includegraphics[width=4.7 in,height=2.85in]{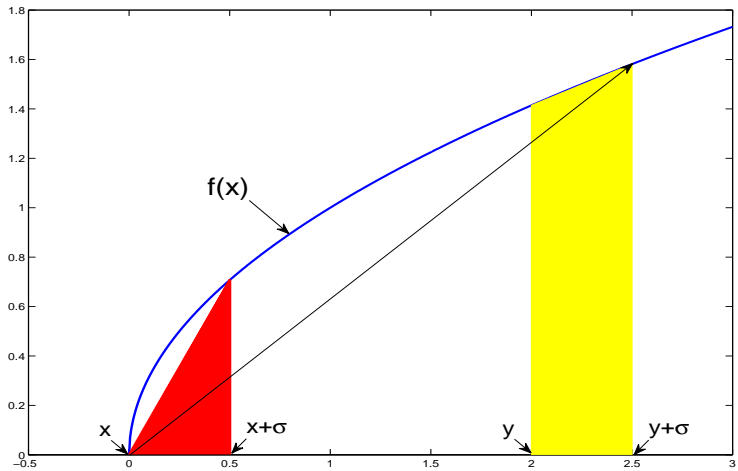}
\caption{\label{fig:1} graph of $f(x)=\sqrt{x}$ when $\sigma=0.5$}
\end{center}
\end{figure}

This picture gives us the idea of Lemma 1 that the slope of an increasing concave function of two points with the same distance of the $x$ coordinate will decrease.
In addition, the monotonicity of the change of slopes will still be valid for other monotone convex or concave functions.
\begin{lemma}\label{l1}
 If $f$ is monotone and concave or convex on an interval $I$ of $\m{R}$, then $G_{\sigma}(x)=|f(x+\sigma)-f(x)|$ is monotone with respect to $x$, where $\sigma$ is any positive real number.
\end{lemma}
\pf We only demonstrate the case that $f(x)$ is monotone increasing and concave, other situations can be proved similarly. Obviously, in this case, $G_{\sigma}(x)=f(x+\sigma)-f(x)$ for $\sigma>0$.
We will show that $G_{\sigma}(x)$ is decreasing on $I$ with respect to $x$.

Suppose $x,y+\sigma \in I$, and $x<y$. Since $f(x)$ is concave, for any $\theta\in(0,1)$, we get
\beq{1l1}
 f(\theta a+(1-\theta)b)\geq \theta f(a)+(1-\theta)f(b), \,\,\,\,\forall a,b\in I.
\eeq

1) Set $a=x$, $b=y+\sigma$, $\theta=\frac{y-x}{y-x+\sigma}$: By \bref{1l1}, we have
$$
f(x+\sigma)\geq \frac{y-x}{y-x+\sigma} f(x) +\frac{\sigma}{y-x+\sigma} f(y+\sigma),
$$
which is equivalent to
\beq{1l2}
\frac{f(y+\sigma)-f(x)}{y-x+\sigma} \leq \frac{f(x+\sigma)-f(x)}{\sigma}.
\eeq

2) Set $a=x$, $b=y+\sigma$, $\theta=\frac{\sigma}{y-x+\sigma}$: Similarly, we can get 
\beq{1l3}
 \frac{f(y+\sigma)-f(y)}{\sigma}\leq \frac{f(y+\sigma)-f(x)}{y-x+\sigma}.
\eeq
By \bref{1l2} and \bref{1l3}, we obtain
$$
 \frac{f(y+\sigma)-f(y)}{\sigma}\leq  \frac{f(x+\sigma)-f(x)}{\sigma}.
$$
Thus, $G_{\sigma}(y) \leq G_{\sigma}(x)$, which implies that $G_{\sigma}(x)$ is a monotone decreasing function.

\rem Lemma \ref{l1} is one special case of a classical property of convex or concave functions on an interval of $\m{R}$, see \cite{convex}.

\begin{lemma}\label{l2}
For a monotone concave or convex function $f(x)$ defined on an interval $I$ of $\m{R}$, if $f$ is uniformly continuous, then $f$ is absolutely continuous.
\end{lemma}
\pf Since $f$ is uniformly continuous on $I$, for any $\epsilon>0$, there exists $\delta>0$
such that for any $x,y\in I$ and $|x-y|<\delta$, we have $$|f(x)-f(y)|<\epsilon.$$
Now, we consider every finite collection $\{(x_i, y_i)\}_{i=1}^{n}$ of nonoverlapping subintervals of $I$ with $x_i<x_{i+1}$, and
$$\sum_{i=1}^{n} (y_i-x_i)<\delta.$$
Define $\sigma_i=y_i-x_i>0$, then $$|f(y_i)-f(x_i)|=G_{\sigma_i}(x_i).$$
By Lemma \ref{l1}, we know that $G_{\sigma}(x)$ is monotone.

1) $G_{\sigma}(x)$ is monotone deceasing.
Since $(x_2,y_2)$ and $(x_1,y_1)$ are nonoverlapping, we get $$G_{\sigma_2}(x_2)\leq G_{\sigma_2}(y_1)=G_{\sigma_2}(x_1+\sigma_1) .$$
Inductively, for $i\geq2$, we will obtain $$G_{\sigma_i}(x_i)\leq G_{\sigma_i}(x_1+\sum_{j=1}^{i-1}\sigma_j). $$
Therefore, we have
\bea
\sum_{i=1}^{n}|f(y_i)-f(x_i)|&=&\sum_{i=1}^{n} G_{\sigma_i}(x_i)\leq G_{\sigma_1}(x_1)+\sum_{i=2}^{n} G_{\sigma_i}(x_1+\sum_{j=1}^{i-1}\sigma_j) .
\eea
Define $z_1=x_1$, and $z_i=x_1+\sum_{j=1}^{i-1}\sigma_j$ $(2\leq i\leq n+1)$. Since $f$ is monotone, the above inequality is equivalent to
$$ \sum_{i=1}^{n}|f(y_i)-f(x_i)|\leq \sum_{i=1}^{n} |f(z_{i+1})-f(z_i)|=|f(z_{n+1})-f(z_{1})|.$$
In addition, $|z_{n+1}-z_{1}|=\sum_{i=1}^{n} \sigma_i =\sum_{i=1}^{n} (y_i-x_i)<\delta $, since $f$ is uniformly continuous, we obtain
$$\sum_{i=1}^{n}|f(y_i)-f(x_i)|\leq |f(z_{n+1})-f(z_{1})|<\epsilon .$$ Hence, $f$ is also absolutely continuous on $I$.

2) $G_{\sigma}(x)$ is monotone increasing. The strategy is quite
similar to the previous one; however, we just fix $(x_n,y_n)$ first
and define $z_{n+1}=y_n$, $z_i=y_n-\sum_{j=i}^{n} \sigma_{j}$ for
$1\leq i\leq n$. Similarly, we have
$$\sum_{i=1}^{n}|f(y_i)-f(x_i)|\leq \sum_{i=1}^{n}|f(z_{i+1})-f(z_i)|=|f(z_{n+1})-f(z_{1})|<\epsilon .$$
Therefore, $f$ is absolutely continuous on $I$.


\rem  The idea of this lemma is to glue disjoint subintervals
$(x_i,y_i)$ together as one subinterval then apply the property of
uniform continuity. By Lemma \ref{l2}, it is clear that
$f(x)=\sqrt{x}$ is absolutely continuous on $[0,c]$. In addition, we
can consider more general functions that oscillate \,\,finite times,
and on each monotone subinterval, they also admit convexity. If
these functions are uniformly continuous, then they are also
absolutely continuous by utilizing the same strategy in the proof of
Lemma \ref{l2} on each subinterval.

\medskip\noindent{\it Proof of Theorem 1.} If $f$ is not monotone on
$[a_i,a_{i+1}]$, then we can split $[a_i,a_{i+1}]$ into two
subintervals such that on each of them, $f$ is monotone, since $f$
assumes convexity. If we relabel them, then on each $[a_i,a_{i+1}]$,
$f$ is monotone and concave or convex. The following proof is based
on this situation.

Since $f$ is uniformly continuous on $I$, for any $\epsilon/N>0$,
there exists $\delta>0$
such that for any $x,y\in I$ and $|x-y|<\delta$, we have $$|f(x)-f(y)|<\epsilon/N.$$ 
Choose $\delta_1<\min\{a_{1}-a_{0},\dots,a_{N}-a_{N-1}, \delta\}$, then for every finite collection $\{(x_i, y_i)\}_{i=1}^{n}$
of nonoverlapping subintervals of $I$ with $x_i<x_{i+1}$ and $\sum_{i=1}^{n} (y_i-x_i)<\delta_1,$ each subinterval $(x_i, y_i)$
can contain at most one $a_j$. For such interval  $(x_i, y_i)$ containing $a_j$, by triangle inequality,  we have
\beq{t1}
|f(y_i)-f(x_i)|\leq |f(y_i)-f(a_j)|+|f(a_j)-f(x_i)|.
\eeq
Now, we still treat $(x_i, a_j)$ and $(a_j, y_i)$ as two ``nonoverlapping subintervals''. Then, we relabel all the subintervals as $\{(x_j, y_j)\}_{j=1}^{m}$ ($n\leq m\leq 2n$).
Through the above strategy, each new $(x_j, y_j)$ will lie exactly in one $[a_i,a_{i+1}]$ and
$$
\sum_{j=1}^{m} (y_j-x_j)=\sum_{i=1}^{n} (y_i-x_i)<\delta_1.
$$

For all the new subintervals that lie in  $[a_i,a_{i+1}]$ {\small$(i=0,\dots,N-1)$},
 since $f$ is also monotone and concave or convex and
$$ \sum_{(x_j,y_j)\subset [a_i,a_{i+1}]} (y_j-x_j)<\delta_1<\delta, $$ by the method in Lemma \ref{l2}, we obtain that
$$
\sum_{(x_j,y_j)\subset [a_i,a_{i+1}]} |f(y_j)-f(x_j)|<\epsilon/N.
$$
Therefore,
\beq{t2}
\sum_{j=1}^{m}|f(y_j)-f(x_j)|=\sum_{i=0}^{N-1} {\sum_{(x_j,y_j)\subset [a_i,a_{i+1}]} |f(y_j)-f(x_j)|}<\epsilon. 
\eeq
By \bref{t1} and \bref{t2}, we get
$$
\sum_{i=1}^{n}|f(y_i)-f(x_i)|\leq \sum_{j=1}^{m}|f(y_j)-f(x_j)|<\epsilon.
$$
Hence, $f$ is also absolutely continuous on $I_{a,b}$.

\bigskip
\noindent {\bf Question:} We consider the case on the real line;
however, one can think about the situation for the multidimensional
Euclidean space.



\bigskip

\noindent\small{\it Department of Mathematics and Statistics,
Memorial University of Newfoundland, NL, Canada.}


\noindent\small{\it kyang@mun.ca}

\medskip\medskip
\noindent\small{\it Applied Mathematical and Computational Sciences,
University of Iowa, IA, USA.}

\noindent\small{\it chenhzhu@math.uiowa.edu}

\end{document}